\documentstyle[11pt,leqno,amscd,amssymb,pstricks,latexsym,amsbsy,xypic,mathrsfs,verbatim]{amsart}

\newcommand{\scr}[1]{\mathscr #1}
\newtheorem{thm}{Theorem}[section]
\newtheorem{lem}[thm]{Lemma}
\newtheorem{cor}[thm]{Corollary}
\newtheorem{prop}[thm]{Proposition}

 \theoremstyle{definition}

\newtheorem{defn}[thm]{Definition}

\newtheorem{rems}[thm]{Remark}

\numberwithin{equation}{thm}

\def\s{{S_\vartriangle(2,2)}}
\def\fs{{\bar{S}_\vartriangle(2,2)}}
\def\en{\bar{e}_\nu}

\def\t1{{S_\vartriangle(2,r)}}

\def\sn{{S_\vartriangle(n,r)}}
\def\bbz{{\mathbb{Z}}}\def\bbq{{\mathbb{Q}}}
\def\bbn{{\mathbb{N}}}
\def\lb{{\lambda}}
\def\c2{{\Theta_\vartriangle(2,2)}}
\def\fkS{{\frak S}}
\def\fsr{{\frak S_{\vartriangle}}}

\def\fii{{\frak S_{\vartriangle,\underline{i}}}}

\def\fiii{{\frak S_{\underline{i}}}}
\def\fjf{{\frak S_{\underline{j}}}}

\def\fj{{\frak S_{\vartriangle,\underline{j}}}}
\def\fjl{{\frak S_{\vartriangle,\underline{j},\underline{l}}}}

\def\fjlf{{\frak S_{\underline{j},\underline{l},\varepsilon'}}}

\def\fijf{{\frak S_{\underline{i},\underline{j},\varepsilon}}}
\def\fijff{{\frak S_{\underline{i},\underline{j}}}}

\def\fij{{\frak S_{\vartriangle,\underline{i},\underline{j}}}}

\def\fil{{\frak S_{\vartriangle, \underline{i}, \underline{l}\delta}}}

\def\fijl{{\frak S_{\vartriangle, \underline{i},\underline{j}, \underline{l}\delta}}}

\def\ff{{\frak S_{\vartriangle, \underline{i}\delta, \underline{l}}}}
\def\fg{{\frak S_{\vartriangle, \underline{i}\delta,\underline{j}, \underline{l}}}}

\def\ffj{{\frak S_{\underline{j}}}}
\def\ffjl{{\frak S_{\underline{j},\underline{l},\varepsilon'}}}
\def\fflj{{\frak S_{\underline{i},\underline{j},\varepsilon}}}

\def\ww{{\frak S_{\underline{i},\underline{l}\delta,\varepsilon'\delta+\varepsilon}}}
\def\vv{{\frak S_{\underline{i},\underline{j},\underline{l}\delta,\varepsilon'\delta,\varepsilon}}}

\def\wv{{\frak S_{\underline{i}\delta,\underline{l},\varepsilon'+\varepsilon\delta}}}
\def\vw{{\frak S_{\underline{i}\delta,\underline{j},\underline{l},\varepsilon',\varepsilon\delta}}}

\def\xij{{ \xi_{\underline{i}, \underline{j}+n\varepsilon}}}
\def\xjl{{ \xi_{\underline{j}, \underline{l}+n\varepsilon'}}}

\def\xil{{ \xi_{\underline{i}, \underline{l}\delta+n(\varepsilon'\delta+\varepsilon)}}}

\def\yil{{ \xi_{\underline{i}\delta, \underline{l}+n(\varepsilon'+\varepsilon\delta)}}}

\def\c{{e_C}}

\def\lr{{\Lambda_\vartriangle(n,r)}}

\def\diag{{\rm diag}}

\def\h{H_0(r)}

\def\bbf{{\mathbb F}}
\def\wfkF{\fkF_\vartriangle}

\def\fkF{{\frak F}}\def\qq{{\boldsymbol q}}
\def\udim{{\underline {\rm dim}}\,}

\def\wfkFn{\fkF_{\vartriangle, n}}

\def\cO{{\cal O}}

\def\row{{\rm row}} \def\col{{\rm col}}

\def\cnr{{\Theta_\vartriangle(n,r)}}
\def\pcnr{{\Theta_\vartriangle^+(n,r)}}
\def\ncnr{{\Theta_\vartriangle^-(n,r)}}

\def\ss{{\widetilde{S}_\vartriangle(n,r)
}}

\def\gp#1#2{\left({#1\atop #2}\right)}
\def\ggp#1#2{\left[\kern-3.2pt\left[{#1\atop #2}\right]\kern-3.2pt\right]}

\def\uj{\underline{j}}

\def\ui{\underline{i}}\def\ul{\underline{l}}
\def\um{\underline{m}}

\def\c{\mathcal{C}}

\def\qq{{\boldsymbol q}}

\def\hr{H_\vartriangle(r)}
\def\h2{H_\vartriangle(2)}

\begin{document}

\title{Affine cellularity of $\s$}

\author{Guiyu Yang}

\thanks{The author is supported by the Natural Science Foundation of
China (Grant no. 11201269 and 11126122),  Natural Science
Foundation of Shandong Province (Grant no. ZR2011AQ004), and
  Young Scholars Research Fund of Shandong University of
Technology.}

\address{School of Science, Shandong University of Technology, Zibo
255049, China} \email{yanggy@@mail.bnu.edu.cn}




\begin{abstract}
In this paper we prove that the affine Schur algebra $\s$ is
affine cellular over $\bbq$. As an application, we prove it is of
finite global dimension.

\end{abstract}

\maketitle

\section{Introduction}
Affine Schur algebras have several equivalent definitions given by
\cite{GV}\cite{L}\cite{rmg}\cite{vv} and \cite{Ya1} respectively.
There are clear correspondences between these definitions on basis
elements. Affine Schur algebras play a central role in linking the
representations of affine quantum groups and affine Hecke
algebras.

Affine cellular algebras are introduced by Koenig and Xi in
\cite{kx}. They extend the framework of cellular algebras due to
Graham and Lehrer to affine cellular algebras which are not
necessarily finite dimensional over a field. Many examples such as
affine Temperley-Lieb algebras and affine Hecke algebras of type
$A$ when the parameter $q$ is not a root of poincare polynomial
are proved to be affine cellular in \cite{kx}. Recently A. S.
Kleshchev, J. W. Loubert, V. Miemietz and J. Guilhot prove that
KLR algebras of type $A$ and affine Hecke algebras of rank two are
affine cellular in \cite{klm} and \cite{gm} respectively.

The aim of this paper is to prove that the affine Schur algebra
$\s$ in case $q=1$ is an affine cellular algebra over $\bbq$. We
use the equivalent definitions of $\sn$ given by \cite{GV}\cite{L}
and \cite{Ya1} respectively. By using the multiplication formulas
given in \cite{L} \cite{DDF} and \cite{Ya}, we investigate the
ideal generated by the idempotent corresponding to a particular
diagonal matrix and construct a chain of idempotent ideals in
$\s$. Then we prove that this chain affords an affine cellular
structure for $\s$. As an application, we prove that $\s$ over
$\bbq$ is of finite global dimension.

The paper is organized as follows. We recall equivalent
definitions of affine Schur algebras and some multiplication
formulas in section 2. In section 3 we introduce affine cellular
algebras defined by Koenig and Xi. In section 4 we prove that $\s$
has an affine cellular structure and is of finite global
dimension.

\section{The affine Schur algebra}

In this section we introduce equivalent definitions of affine
Schur algebras given by \cite{GV} \cite{L} and \cite{Ya1}. We give
correspondences of the basis elements between these definitions.
At the end of this section, we recall some multiplication formulas
given by \cite{DDF} \cite{DF} \cite{L} and \cite{Ya1}
respectively.

 First we give a geometric
definition of affine quantum Schur algebras introduced by
Ginzburg--Vasserot \cite{GV} and Lusztig\cite{L}.
 Let $\bbf$ be a field and let
$\bbf[x,x^{-1}]$ be the Laurent polynomial ring in indeterminate
$x$. Fix an $\bbf[x,x^{-1}]$-free module $V$ of rank $r\geq 1$. A
lattice in $V$ is a free $\bbf[x]$-submodule $L$ of $V$ such that
$V=L\otimes_{\bbf[x]}\bbf[x,x^{-1}]$.

 Let $\wfkF=\wfkFn$ be the set of all cyclic flags $L=(L_i)_{i\in\bbz}$
of lattices, where each $L_i$ is a lattice in $V$ such that
$L_{i-1}\subseteq L_i$ and $L_{i-n}=x L_i$ for all $i\in \bbz.$
The group $G$ of automorphisms of the the $\bbf[x,x^{-1}]$-module
$V$ acts on $\wfkF$ by $g\cdot L=(g(L_i))_{i\in\bbz}$ for $g\in G$
and $L\in \wfkF$. Thus, the map
$$\phi:\wfkF\rightarrow \lr,\;\; L\rightarrow
\udim L=(\dim_\bbf L_i/L_{i-1})_{i\in \bbz}$$ induces a bijection
between the set of $G$-orbits in $\wfkF$ and $\lr$, where
$$\lr:=\{(\lb_i)_{i\in\bbz}\mid \lb_i\in\bbn,
\sum_{i=1}^n\lb_i=r \;\mbox{and}\;\lb_i=\lb_{i-n}\; \mbox{for}\;
i\in \bbz\}.$$

Let $$\Lambda(n,r):=\{(\lb_1,\ldots,\lb_n)~|~\lb_i\in\bbn,~
\sum_{1\leq i\leq n}\lb_i=r\}.$$  We usually identify
$\Lambda(n,r)$ with $\lr$ via the following bijection:
$$b: \lr\longrightarrow \Lambda(n,r),~~\lb\longrightarrow (\lb_1,\ldots,\lb_n).$$

 The group $G$ also acts diagonally on $\wfkF\times \wfkF$ by $g(L,
L')=(gL, gL')$, where $g\in G$ and $L, L'\in \wfkF$. By
\cite[1.5]{L}, there is a bijection between the set of $G$-orbits
in $\wfkF\times\wfkF$ and the set $\cnr$ by sending $(L,L')$ to
$A=(a_{i,j})_{i,j\in\bbz}$, where
\begin{equation}\label{aij'} a_{i,j}=\dim_\bbf\frac{L_i\cap
L_j'}{L_{i-1}\cap L_j'+L_i\cap L_{j-1}'}\;\;\text{for $
i,j\in\bbz$},
\end{equation}

\begin{equation}\label{ba}\Theta_\vartriangle(n,r):=\{A=(a_{i,j})_{i,j\in\bbz}\in
M_{\vartriangle,n}(\bbn)\mid  \sum_{1\leq i\leq
n\atop{j\in\bbz}}a_{i,j}=\sum_{1\leq j\leq
n\atop{i\in\bbz}}a_{i,j}=r \}\end{equation} and
$M_{\vartriangle,n}(\bbn)$ is the set of all $\bbz\times\bbz$
 matrices $A=(a_{i,j})_{i,j\in\bbz}$ with $a_{i,j}\in\bbn$ such that
\begin{enumerate}
\item[(a)] $a_{i,j}=a_{i+n,j+n} \;\mbox{for}\; i, j \in \bbz\;; $

\item[(b)]$\mbox{for every}\; i\in \bbz,\;\mbox{the
set}\;\{j\in\bbz\mid a_{i,j}\neq 0\}\;\mbox{is finite}\;. $
\end{enumerate}
Let $\cO_A$ denote the orbit in $\wfkF\times\wfkF$ corresponding
to $A$. If $(L,L')\in\cO_A$, then $\row(A)=\udim L$ and
$\col(A)=\udim L'$, where $\row(A)=(\sum_{j\in\bbz}
a_{i,j})_{i\in\bbz}$ and $\col(A)=(\sum_{i\in\bbz}
a_{i,j})_{j\in\bbz}$.

Assume now that $\bbf=\bbf_q$ is a finite field of $q$ elements
and write $\wfkF(q)$ for $\wfkF$. For any fixed
$(L,L'')\in\cO_{A''}$, let $c_{A,A',A'';q}=|\{L'\in\wfkF(q)~|~
(L,L')\in\cO_A, (L',L'')\in\cO_{A'}\}|$. Clearly, $c_{A,A',A'';q}$
is independent of the choice of $(L,L'')$, and a necessary
condition for $c_{A,A',A'';q}\neq 0$ is that
\begin{equation}\label{ee}\col(A)=\row(A'),~ \row(A)=\row(A'')
~\text{and}~\col(A')=\col(A'').\end{equation}

Let ${\scr A}=\bbz[\qq]$ be the polynomial ring with indeterminate
$\qq$. By \cite{L}, there is a polynomial $p_{A,A',A''}\in
\bbz[\qq]$ in $\qq$ such that for each finite field $\bbf$ with
$q$ elements, $c_{A,A',A'';q}=p_{A,A',A''}(q)$.

\begin{defn} [{\cite{L}},\cite{GV}] \label{def-q-Schur}
The affine quantum Schur algebra $\sn$ is the free
$\bbz[\qq]$-module with basis $\{e_A~|~ A\in\cnr\}$, and
multiplication defined by $$e_{A}\cdot e_{A'}=\begin{cases}\sum_{A''\in\cnr}p_{A,A',A''}e_{A''},~~\text{if}~ \col(A)=\row(A'),\\
              0,~~\text{otherwise}.
\end{cases} $$
\end{defn}

As in the finite case, for each $\lb\in\lr$, define
\diag$(\lb)=(\delta_{i,j}\lb_i)_{i,j\in\bbz}\in
\Theta_\vartriangle(n,r)$, and $e_\lb=e_{\text{\diag}(\lb)}$. It
is easy to see that for each $A\in\cnr$,
\begin{equation}\label{idempotents}
e_{\lb}e_A=\left\{\begin{array}{ll}
                      e_A,\;&\text{if $\lb=\row(A)$};\\
                     0,
                     &\text{otherwise}\end{array}\right.\;\;\text{and}\;\;
e_Ae_{\lb}=\left\{\begin{array}{ll}
                      e_A,\;&\text{if $\lb=\col(A)$};\\
                     0, &\text{otherwise.}\end{array}\right.
\end{equation}
 Thus, $\sum_{\lb\in\lr}e_{\lb}$ is the unity
 of $\sn$. By specializing $v$ to $1$ in definition \eqref{def-q-Schur}, we
get the affine Schur algebra $\sn$ over $\bbz$, which we still
denote by $\sn$.

Now we introduce an algebraic definition of affine Schur algebras
given by \cite{Ya1}. Let $\fkS=\fkS_r$ denote the symmetric group
on $r$ letters and let $\fsr=\fkS\ltimes \bbz^r$ denote the
extended affine Weyl group of type $A_{r-1}$. For a set $S$, we
denote by $I(S,r)$ the set
$\{\underline{i}=(i_1,\ldots,i_r)~|~i_t\in S, t=1,2,\ldots,r\}$.
We will denote the set $I(S,r)$ by $I(n,r)$ if
$S=\{1,2,\ldots,n\}$. Now $\fkS$ acts on $I(n,r)$ by place
permutation. $\fsr$ acts on $I(\bbz,r)$ on the right with $\fkS$
acting by place permutation and $\bbz^r$ acting by shifting, i.e.
$\underline{i}(\sigma, \varepsilon)=\underline{i}+n\varepsilon$
for $\underline{i}\in I(\bbz,r)$, $\sigma\in\fkS$ and
$\varepsilon\in \bbz^r$. Note that this action depends on the
number $n$ and $\fsr$ acts diagonally on $I(\bbz,r)\times
I(\bbz,r)$. For $(\underline{i}, \underline{j})\in I(\bbz,r)\times
I(\bbz,r)$ and $(\underline{k}, \underline{l})\in I(\bbz,r)\times
I(\bbz,r)$, we identify $\xi_{\underline{i}, \underline{j}}$ and
$\xi_{\underline{k}, \underline{l}}$ if and only if
$(\underline{i}, \underline{j})\sim _\fsr(\underline{k},
\underline{l})$, i.e $(\underline{i}, \underline{j})$ and
$(\underline{k}, \underline{l})$ are in the same orbit.

\begin{defn} [{\cite{Ya1}}] \label{def-q-Schur2}
The affine Schur algebra $\ss$ is defined to be the $\bbz$-algebra
with basis $\{\xi_{\underline{i}, \underline{j}}~|~\underline{i},
\underline{j} \in I(\bbz, r)\}$ and multiplication given by the
following rule:

$$\xi_{\underline{i},
\underline{j}}\xi_{\underline{k},
\underline{l}}=\sum_{(\underline{p}, \underline{q})\in
I(\bbz,r)\times I(\bbz,r)/\fsr }C(\underline{i},
\underline{j},\underline{k}, \underline{l},\underline{p},
\underline{q})\xi_{\underline{p}, \underline{q}},$$ where
$C(\underline{i}, \underline{j},\underline{k},
\underline{l},\underline{p}, \underline{q})=|\{\underline{s}\in
I(\bbz, r)~|~(\underline{i}, \underline{j})\sim
_\fsr(\underline{p}, \underline{s}),~(\underline{s},
\underline{q})\sim _\fsr(\underline{k}, \underline{l})\}|.$

\end{defn}

\begin{rems}\label{r1}
 There is an $\bbz$-algebra isomorphism between $\sn$ and
$\ss$. The correspondence $\varphi$ between the basis elements is
given by
$$\varphi(\xi_{\underline{i},
\underline{j}})=e_A,$$ where $A=(a_{x,y})_{x,y\in\bbz}\in\cnr$,
$\ui, \uj\in I(\bbz, r)$ and
$$a_{x,y}=|\{s~|~i_s=x, j_s=y, 1\leq s\leq r\}|.$$
\end{rems}

The $\bbz$-algebra isomorphism above can be generalized to a
$\bbq$-algebra isomorphism. In this paper, we are mainly
concentrated in the affine Schur algebra $\s$ over $\bbq$. For
simplicity, we usually identify the basis elements occurring in
the above equivalent definitions. And we always consider $\s$ as
an algebra over $\bbq$.

\begin{rems}\label{rd}
By section 5.2 in \cite{DDF}, $\sn$ is a graded algebra over
$\bbq$, i.e.
$$\sn=\bigoplus_{m\in\bbz}\sn_m.$$
Let $\pcnr=\{A=(a_{i,j})\in\cnr~|~a_{i,j}=0 ~\text{for}~ i>j\}$
and $\ncnr=\{A=(a_{i,j})\in\cnr~|~a_{i,j}=0 ~\text{for}~ i<j\}$.
Then the degree gr$(e_A)$ for $A\in\pcnr$ is defined by
$$\text{gr}(e_A)=\sum_{i<j, 1\leq i\leq n}a_{i,j}(j-i),$$
and the degree gr$(e_A)$ for $A\in\ncnr$ is defined by
$$\text{gr}(e_A)=\sum_{i>j, 1\leq i\leq n}-a_{i,j}(i-j).$$

\end{rems}

Now we give some useful multiplication formulas which are given by
\cite{L}, \cite{DDF}, \cite{DF} and \cite{Ya1} respectively.

For $i, j\in \bbz$, let $E_{i,j}^\vartriangle\in
M_{\vartriangle,n}(\bbn)$ be the matrix $(e_{k,l}^{i,j})_{k,l\in
\bbz}$ defined by

$$e_{k,l}^{i,j}=\begin{cases}1,~~\text{if}~ k=i+sn, l=j+sn ~\text{for some}~ s\in\bbz,\\
              0,~~\text{otherwise}.
\end{cases}$$

\begin{prop}[\cite{L},\cite{DDF},\cite{DF}]\label{mu1}
Let $1\leq h\leq n, A\in\cnr$ and $\lb=\row(A)$. Let
$B_m=\diag(\lb)+mE_{h,h+1}^\vartriangle-mE_{h+1,h+1}^\vartriangle$
and
$C_m=\diag(\lb)-mE_{h,h}^\vartriangle+mE_{h+1,h}^\vartriangle$.
Then in $\sn$

$$(1) ~e_{B_m}\cdot e_A=\sum_{t\in \Lambda(\infty, m)\atop
{\forall u\in\bbz, t_u\leq a_{h+1,u}}}\prod_{u\in \bbz}\gp
{a_{h,u}+t_u} {t_u}e_{A+\sum_{u\in\bbz}t_u
(E_{h,u}^\vartriangle-E_{h+1,u}^\vartriangle)},$$

$$(2) ~e_{C_m}\cdot e_A=\sum_{t\in \Lambda(\infty,
m)\atop {\forall u\in\bbz, t_u\leq a_{h,u}}}\prod_{u\in \bbz}\gp
{a_{h+1,u}+t_u} {t_u}e_{A-\sum_{u\in\bbz}t_u
(E_{h,u}^\vartriangle-E_{h+1,u}^\vartriangle)},$$ where
$\Lambda(\infty, m)=\{\lb=(\lb_i)_{i\in\bbz}~|~\lb_i\in\bbn ,
\sum_i\lb_i=m\}$.

\end{prop}

Symmetrically, we can get the following multiplication formulas.
\begin{prop}\label{mu11}
Let $1\leq h\leq n, A\in\cnr$ and $\lb=\col(A)$. Let
$B_m=\diag(\lb)+mE_{h,h+1}^\vartriangle-mE_{h,h}^\vartriangle$ and
$C_m=\diag(\lb)-mE_{h+1,h+1}^\vartriangle+mE_{h+1,h}^\vartriangle$.
Then in $\sn$

$$(1) ~ e_A\cdot e_{B_m}=\sum_{t\in \Lambda(\infty, m)\atop
{\forall u\in\bbz, t_u\leq a_{u,h}}}\prod_{u\in \bbz}\gp
{a_{u,h+1}+t_u} {t_u}e_{A+\sum_{u\in\bbz}t_u
(E_{u,h+1}^\vartriangle-E_{u,h}^\vartriangle)},$$

$$(2) ~e_A \cdot e_{C_m}=\sum_{t\in \Lambda(\infty,
m)\atop {\forall u\in\bbz, t_u\leq a_{u,h+1}}}\prod_{u\in \bbz}\gp
{a_{u,h}+t_u} {t_u}e_{A-\sum_{u\in\bbz}t_u
(E_{u,h+1}^\vartriangle-E_{u,h}^\vartriangle)},$$ where
$\Lambda(\infty, m)=\{\lb=(\lb_i)_{i\in\bbz}~|~\lb_i\in\bbn ,
\sum_i\lb_i=m\}$.

\end{prop}

\begin{prop}[\cite{DDF}]\label{mu3}

Let $1\leq h\leq n, A\in\cnr$ and $\lb=\row(A)$. Let
$D_m=\diag(\lb)-E_{h,h}^\vartriangle+E_{h,h+mn}^\vartriangle$.
Then in $\sn$
$$(1) ~e_{D_m}\cdot e_A=\sum_{u\in\bbz\atop
{ a_{h,u}\geq 1}}(a_{h,u+mn}+1)e_{A+
(E_{h,u+mn}^\vartriangle-E_{h,u}^\vartriangle)},~~~~~~~~~~~~~~~~~~~~~~$$
where $m\in\bbz\setminus\{0\}$.

\end{prop}

\begin{prop}[\cite{Ya1}]\label{mu2} Let $\ui, \uj, \underline{k}, \ul\in I(\bbz,
r)$. We have the following equations in $\sn$.

 $(1)
~~\xi_{\ui, \uj }\xi_{\underline{k}, \ul}=0$ unless $\uj\sim_\fsr
\underline{k}$.

 $(2)
~~\xi_{\ui, \ui }\xi_{\underline{i}, \uj}=\xi_{\underline{i},
\uj}=\xi_{\ui, \uj }\xi_{\underline{j}, \uj}$.

  $(3)
~~\sum_{\ui\in I(n,r)/\fkS}\xi_{\ui, \ui }$ is a decomposition of
unity into orthogonal idempotents.

\end{prop}

\begin{prop}[\cite{Ya1}, \cite{Ya}]\label{mu4}

$$\aligned(1)~~~ \xi_{\underline{i},
\underline{j}}\xi_{\underline{j},
\underline{l}}=\sum_{\delta\in\fjl\setminus \fj/ \fij
}\biggl[\fil:\fijl\biggr]\xi_{\underline{i},
\underline{l}\delta}\\
=\sum_{\delta\in\fij\setminus \fj/ \fjl
}\biggl[\ff:\fg\biggr]\xi_{\underline{i}\delta, \underline{l}},\\
\endaligned$$
 where $\ui, \uj, \ul$ are in $I(\bbz,r)$,
$\fii$ is the stabilizer subgroup of $\ui$ in $\fsr$ and $\fij$ is
the stabilizer of $\ui$ and $\uj$ in $\fsr$, i.e. $\fij=\fii\cap
\fj$, etc, $\fjl\setminus \fj/ \fij$ denotes a representative set
of double cosets.

$$\aligned~~~(2) \xij\xjl=\sum_{\delta\in\ffjl\setminus \ffj/ \fflj
}\biggl[\ww:\vv\biggr]\xil\\
 =\sum_{\delta\in\fflj\setminus \ffj/ \ffjl
}\biggl[\wv:\vw\biggr]\yil\\
\endaligned$$
 where $\ui, \uj, \ul$ are in $I(n,r)$, $\varepsilon,
 \varepsilon'\in \bbz^r$,
$\fiii$ is the stabilizer subgroup of $\ui$ in $\fkS$ and $\fijff$
is the stabilizer of $\ui$ and $\uj$ in $\fkS$, i.e.
$\fijff=\fiii\cap \fjf$, etc, $\fjlf\setminus \fjf/ \fijf$ denotes
a representative set of double cosets.

\end{prop}
\section{The affine cellular algebra}

In this section let $k$ be a noetherian domain. A commutative
$k$-algebra $B$ is called affine if it is a quotient of a
polynomial ring $k[x_1,\ldots,x_t]$ in finitely many variables. A
$k$-linear anti-automorphism $\tau$ of a $k$-algebra $A$ with
$\tau^2=id_A$ will be called a $k$-involution on $A$.

\begin{defn} [{\cite{kx}}] \label{def-aff}
Let $A$ be a unitary $k$-algebra with a $k$-linear involution
$\tau$ on $A$. A two-sided ideal $J$ in $A$ is called an affine
cell ideal if the following conditions are satisfied:

 (1) $\tau(J)=J$.

(2) There is a free $k$-module $V$ of finite rank and an affine
algebra $B$ with a $k$-involution $\sigma$ such that
$\triangle=V\otimes_kB$ is an $A$-$B$-bimodule, where the right
$B$-module structure is induced by that of the regular right
$B$-module $B_B$ .

(3) There is an $A$-$A$-bimodule isomorphism $\alpha:J\rightarrow
\triangle\otimes_B\triangle'$, where $\triangle'=B\otimes_kV$ is a
 $B$-$A$-bimodule with the left $B$-structure induced by $_BB$ and the right $A$-module structure is
given as $(b\otimes v)a=p(\tau(a)(v\otimes b)), $ where $p$ is the
switch map: $$p:~~\triangle\otimes_B\triangle'\longrightarrow
\triangle'\otimes_B\triangle,$$$$ x\otimes y\longrightarrow
y\otimes x, ~~\text{for}~~ x\in\triangle ~~\text{and}~~ y\in
\triangle'.$$

(4) There is the following commutative diagram:

$$\xymatrix{ J \ar[d]^-{\tau}\ar[r]^-{\alpha}& \triangle \otimes_B \triangle'\ar[d]^-{
v_1\otimes b_1\otimes_B b_2\otimes v_2\rightarrow v_2\otimes \sigma(b_2)\otimes_B \sigma(b_1)\otimes v_1}\\
J\ar[r]^-{\alpha} &\triangle \otimes_B \triangle' }$$

The algebra $A$ is called affine cellular if and only if there is
a $k$-module decomposition $A=J_1'\oplus J_2'\oplus\cdots J_n'$
with $\tau(J_j')=J_j'$ for each $j$, and $J_i=\oplus_{1\leq l\leq
i}J_l'$ gives a chain of two-sided ideals of $A$: $0=J_0\subset
J_1\subset J_2\subset\cdots\subset J_n=A$, and for each $1\leq
i\leq n$, $J_i'=J_i/J_{i-1}$ is an affine cell ideal of
$A/J_{i-1}$.

\end{defn}

\begin{rems}\label{a1}  Note that
by the definition of affine cellular algebras, quotient algebras
of commutative polynomial rings in finitely many variables over a
noetherian domain $k$ are always affine cellular, with identity
equals $\tau$, $B$ equals itself and $V$ equals the trivial
$B$-module $k$. In particular, the affine Schur algebra
$S_\vartriangle(1,r)$ over $\bbq$ is isomorphic to $\bbq[x_1,
x_2,\ldots, x_{r-1}, x_r, x_r^{-1}]$. So it is affine cellular. In
the following section we will prove that $\s$ is affine cellular
over $\bbq$.

\end{rems}

\section{ affine cellularity of $\s$}

In this section we give a proof of the affine cellularity of $\s$
and prove it is of finite global dimension. First we recall the
definition of affine Hecke algebra $\hr$ over $\bbq[\qq]$. It has
a set of generators $T_i(1\leq i\leq r), T_\rho^{\pm 1}$ with the
following relations:\begin{equation}\label{re}
\begin{cases} T_i^2=(\qq-1)T_i+\qq, \; \;& \text{for} \;1\leq i\leq r;\\
T_iT_j=T_jT_i, \;\;&\mbox{for $1\leq i,j\leq r$ with $|i-j|>1$};\\
T_iT_{i+1}T_i =T_{i+1}T_iT_{i+1},\;\;&\text{for $1\leq i\leq r$, with $|i-j|=1$ and $r\geq 3$};\\
T_\rho T_\rho^{-1}=T_\rho^{-1}T_\rho=1, ~~T_\rho
T_iT_\rho^{-1}=T_{i+1},\;\;&\text{for $1\leq i\leq r$}.
\end{cases}\end{equation}
In the relations above, we identify $T_{r+1}$ with $T_1$. By
specializing $\qq$ to $1$, $\h2$ over $\bbq$ has a presentation as
follows:
\begin{equation}\label{rrr}
\begin{cases} T_1^2=1, \; \;T_2^2=1, ~~T_\rho T_\rho^{-1}=T_\rho^{-1}T_\rho=1;\\
 ~~T_\rho
T_1T_\rho^{-1}=T_{2},~~T_\rho T_2T_\rho^{-1}=T_{1}.
\end{cases}\end{equation}

Let $\lb=(2,0), \mu=(0,2), \nu=(1,1)\in \Lambda(2,2)$. Then
$e_\lb+e_\mu+e_\nu$ is the unity in $\s$. We can identify $\h2$
with the subalgebra $e_\nu\s e_\nu$ of $\s$. The correspondence is
given as follows:

\begin{equation}\label{hs}\varphi: \h2\longrightarrow e_\nu\s
e_\nu,\end{equation}

$$\varphi(1)=e_\nu,~~\varphi(T_1)= e_{(E_{1,2}^\vartriangle+E_{2,1}^\vartriangle)},~~
~~\varphi(T_2)=
e_{(E_{2,3}^\vartriangle+E_{3,2}^\vartriangle)},$$$$
\varphi(T_\rho)=
e_{(E_{1,2}^\vartriangle+E_{2,3}^\vartriangle)},~~
~~\varphi(T_{\rho^{-1}})=
e_{(E_{2,1}^\vartriangle+E_{3,2}^\vartriangle)}.$$

\begin{lem}\label{maa} Let $J$ denote the ideal generated by $e_\lb$ in $\s$, i.e.
$J=\s e_\lb\s$. Then $\fs: =\s/ J$ is isomorphic to
$\bbq[x,x^{-1}]$, where $\bbq[x,x^{-1}]$ is the Laurent polynomial
ring in variable $x$.
\end{lem}

\begin{pf}

By Proposition \eqref{mu1}, we get that
$e_\mu=e_{2E_{2,2}^\vartriangle}=e_{2E_{2,1}^\vartriangle}\cdot
e_{2E_{1,2}^\vartriangle}$. Since $e_{2E_{2,1}^\vartriangle}\cdot
e_\lb=e_{2E_{2,1}^\vartriangle}$ and $e_\lb\cdot
e_{2E_{1,2}^\vartriangle}=e_{2E_{1,2}^\vartriangle}$, we get that
$e_\mu=e_{2E_{2,1}^\vartriangle}\cdot e_\lb\cdot
e_{2E_{1,2}^\vartriangle}\in J$.

 So $\bar{e}_\nu$ is the unity in $\fs$, and $\en
\fs\en=\fs$. Note there is the following isomorphism
$$\en
\fs\en\cong (e_\nu \s e_\nu)/(e_\nu J e_\nu).$$

By identifying $\h2$ with $e_\nu \s e_\nu$ by \eqref{hs}, we prove
that $J$ equals the ideal of $\s$ generated by $T_1+e_\nu$ and
$T_2+e_\nu$.

Let $J'$ be the ideal of $\s$ generated by $T_1+e_\nu$ and
$T_2+e_\nu$. We first show that $J'\subseteq J$. By proposition
\eqref{mu1}, we have
$$e_{(E_{1,1}^\vartriangle+E_{2,1}^\vartriangle)}\cdot
e_\lb=e_{(E_{1,1}^\vartriangle+E_{2,1}^\vartriangle)},~~~~
e_\lb\cdot e_{(E_{1,1}^\vartriangle+E_{1,2}^\vartriangle)}=
e_{(E_{1,1}^\vartriangle+E_{1,2}^\vartriangle)}.$$ Then we get
$$e_{(E_{1,1}^\vartriangle+E_{2,1}^\vartriangle)}\cdot
e_{(E_{1,1}^\vartriangle+E_{1,2}^\vartriangle)}=e_{(E_{1,1}^\vartriangle+E_{2,1}^\vartriangle)}\cdot
e_\lb\cdot e_{(E_{1,1}^\vartriangle+E_{1,2}^\vartriangle)}\in J.$$
Since
$$e_{(E_{1,1}^\vartriangle+E_{2,1}^\vartriangle)}\cdot
e_{(E_{1,1}^\vartriangle+E_{1,2}^\vartriangle)}=
e_{(E_{1,2}^\vartriangle+E_{2,1}^\vartriangle)}+e_\nu=T_1+e_\nu,$$
this proves $T_1+e_{\nu}\in J$. By \eqref{rrr}, $
T_\rho^{-1}(T_1+e_\nu)T_\rho=T_{2}+e_\nu$, then $T_2+e_{\nu}\in
J$. This proves $J'\subseteq J$.

By Proposition \eqref{mu1}, we have
$$e_{(E_{1,1}^\vartriangle+E_{1,2}^\vartriangle)}(T_1+e_\nu)
e_{(E_{1,1}^\vartriangle+E_{2,1}^\vartriangle)}=4e_\lb,$$ then
$e_\lb\in J'$, which proves $J\subseteq J'$.

 Note that $e_\nu J' e_\nu$ equals the ideal of $\h2$ generated by
 $T_1+1$ and $T_2+1$. By $J=J'$ and identifying $\h2$ with $e_\nu
\s e_\nu$, we have  $$(e_\nu \s e_\nu)/(e_\nu J e_\nu)\cong
\h2/(e_\nu J' e_\nu)\cong \bbq[T_\rho,T_\rho^{-1}].$$ This
completes the Lemma.

\end{pf}

\begin{rems}\label{r2}

 By Remarks \eqref{r1}, $T_1=e_{(E_{1,2}^\vartriangle+E_{2,1}^\vartriangle)}$
 corresponds to $\xi_{\ui,\uj}$, where $\ui=(1,2)$ and $\uj=(2,1)$
 are both in $I(2,2)$. Since $(\ui, \uj)\sim_\fsr
 (\uj, \ui)$, by Proposition \eqref{mu4}, $$T_1^2=\xi_{\ui,\uj} \cdot \xi_{\ui,\uj}
 =\xi_{\ui,\uj} \cdot \xi_{\uj,\ui}=\xi_{\ui,\ui}=e_\nu. $$
In the following, we will describe $J$ as a bimodule over $\s$,
and prove it is projective as both left and right $\s$-module.
This affords a solution to bound the global dimension of $\s$.

\end{rems}

The following lemma which is given in \cite{Ya} describes the
 subalgebra $e_\lb \s e_\lb$.

\begin{lem}[\cite{Ya}]\label{m} $e_\lb \s e_\lb$ is isomorphic to the Laurent
polynomial ring $\bbq$$[x_1,x_2,x_2^{-1}]$.
\end{lem}

The isomorphism correspondence between $e_\lb \s e_\lb$ and
$\bbq$$[x_1,x_2,x_2^{-1}]$ is given as follows:

\begin{equation}\label{hs2}\psi:
 e_\lb \s e_\lb\longrightarrow \bbq[x_1,x_2,x_2^{-1}],\end{equation}

$$\psi(e_{(E_{1,1}^\vartriangle+E_{1,3}^\vartriangle)})=x_1,~~\psi(e_{(2E_{1,3}^\vartriangle)})=
x_2,~~ ~~\psi(e_{(2E_{3,1}^\vartriangle)})=x_2^{-1}.$$

\begin{lem}\label{main} Let $B=e_\lb \s e_\lb$. Then $ e_\lb\s$ is a free left
$B$-module of rank four, $\s e_\lb$ is a free right $B$-module of
 rank four.
\end{lem}

\begin{pf}
$ e_\lb\s=\bbq\{e_A~|~
\row(A)=\lb\}=\bbq\{e_{(E_{1,i}^\vartriangle+E_{1,j}^\vartriangle)}~|~i,
j\in\bbz\}$. Identify $B$ with $\bbq[x_1,x_2,x_2^{-1}]$ and by
propositions \eqref{mu3}, the actions of $x_1$ on $ e_\lb\s$ is
given as follows.

$$
x_1\cdot e_{(E_{1,i}^\vartriangle+E_{1,j}^\vartriangle)}=
\begin{cases}
e_{(E_{1,i+2}^\vartriangle+E_{1,j}^\vartriangle)}+
e_{(E_{1,i}^\vartriangle+E_{1,j+2}^\vartriangle)},& i\neq j, ~i\neq j\pm2;\\
e_{(E_{1,i+2}^\vartriangle+E_{1,j}^\vartriangle)},  & i=j ;\\
2e_{(E_{1,i+2}^\vartriangle+E_{1,j}^\vartriangle)}+
e_{(E_{1,i}^\vartriangle+E_{1,j+2}^\vartriangle)},  & i=j-2 ,\\
e_{(E_{1,i+2}^\vartriangle+E_{1,j}^\vartriangle)}+
2e_{(E_{1,i}^\vartriangle+E_{1,j+2}^\vartriangle)},  & i=j+2 .\\
\end{cases}
$$

By Remark \eqref{r1}, identify $x_2$ with $\xi_{\ul,\um}$, where
$\ul=(1,1), \um=(3,3)\in I(\bbz,2)$, and identify
$e_{(E_{1,i}^\vartriangle+E_{1,j}^\vartriangle)}$ with
$\xi_{\underline{l},\underline{u}}$, where $\underline{l}=(1,1),
\underline{u}=(i,j)\in I(\bbz,2)$. By \eqref{mu4} we have
$$
x_2\cdot e_{(E_{1,i}^\vartriangle+E_{1,j}^\vartriangle)}=
e_{(E_{1,i+2}^\vartriangle+E_{1,j+2}^\vartriangle)},$$$$
x_2^{-1}\cdot e_{(E_{1,i}^\vartriangle+E_{1,j}^\vartriangle)}=
e_{(E_{1,i-2}^\vartriangle+E_{1,j-2}^\vartriangle)}.
$$

Now we prove that $ e_\lb\s$ as a left $B$-module is generated by
the set
$$\Pi_1=\{e_{(E_{1,1}^\vartriangle+E_{1,2}^\vartriangle)},~
e_{(E_{1,2}^\vartriangle+E_{1,3}^\vartriangle)},~
e_{(2E_{1,1}^\vartriangle)} ,~ e_{(2E_{1,2}^\vartriangle)} \}.$$

Let $\mathcal{X}$ be the set
$=\{e_{(E_{1,i}^\vartriangle+E_{1,j}^\vartriangle)}~|~i,
j\in\bbz\}$. We prove that each element in $\mathcal{X}$ is
contained in the $B$-linear combinations of elements in $\Pi$.
First if $i=j$, we have $$ e_{(2E_{1,i}^\vartriangle)}=
\begin{cases}
x_2^{k}\cdot e_{(2E_{1,1}^\vartriangle)},  & i=2k+1 ;\\
x_2^{(k-1)}\cdot e_{(2E_{1,2}^\vartriangle)},& i=2k.\\
\end{cases}
$$

If $i\neq j$, we have
$$e_{(E_{1,i}^\vartriangle+E_{1,j}^\vartriangle)}=
\begin{cases}
x_2^k\cdot e_{(E_{1,1}^\vartriangle+E_{1,j-i+1}^\vartriangle)},  & i=2k+1 ;\\
x_2^{(k-1)}\cdot e_{(E_{1,2}^\vartriangle+E_{1,j-i+2}^\vartriangle)},& i=2k.\\
\end{cases}
$$

So we only need to prove that
$X_l=e_{(E_{1,1}^\vartriangle+E_{1,l}^\vartriangle)}$ and
$Y_l=e_{(E_{1,2}^\vartriangle+E_{1,l}^\vartriangle)}$ for
$l\in\bbz $ can be generated by $\Pi_1$.

Note that \begin{equation}\label{h}x_1\cdot X_l=\begin{cases}
X_{l+2}+x_2\cdot X_{l-2},  &l\neq -1, 1, 3 ;\\
X_5+2x_2X_1,& l=3;\\
X_3, &l=1;\\
2X_1+x_2X_{-3}, &l=-1.\\
\end{cases}\end{equation}

 This gives the equalities
\begin{equation}\label{g}X_{l+2}=x_1\cdot X_l-x_2\cdot
X_{l-2}, ~~\text{and}~~ X_{l-2}=x_2^{-1}x_1\cdot X_l-x_2^{-1}\cdot
X_{l+2},\end{equation} for $l\neq -1,3$.
 Since $X_1,
X_2\in\Pi$, by \eqref{h},
$$X_3=x_1X_1,~~X_4=x_1\cdot
X_2-e_{(E_{1,2}^\vartriangle+E_{1,3}^\vartriangle)},$$$$X_5=x_1\cdot
X_3-2x_2\cdot X_1,\;\; X_{-3}=x_2^{-1}x_1\cdot
X_{-1}-2x_2^{-1}\cdot X_1$$ are all generated by $\Pi_1$.

So we can prove each $X_l$ for $l\in\bbz$ is generated by $\Pi_1$
by induction on $l$ for $l\in\bbz$ using \eqref{g}.

Similarly we have
\begin{equation}\label{l}x_1\cdot Y_l=\begin{cases}
Y_{l+2}+x_2\cdot Y_{l-2},  &l\neq 0, 2, 4 ;\\
Y_6+2x_2Y_2,& l=4;\\
Y_4, &l=2;\\
2Y_2+x_2Y_{-2}, &l=0.\\
\end{cases}\end{equation}

 Since $Y_1, Y_2\in\pi$, we can use a similar induction on $l$
 to prove that $Y_l$ for $l\in\bbz$ are all generated by $\Pi_1$.

Now we prove elements in $\Pi_1$ are $e_\lb \s e_\lb$-linearly
independent. By Remark \eqref{rd}, $\s$ is a $\bbz$-graded algebra
over $\bbq$ and the degrees of
$$e_{(E_{1,1}^\vartriangle+E_{1,2}^\vartriangle)},~
e_{(E_{1,2}^\vartriangle+E_{1,3}^\vartriangle)},~
e_{(2E_{1,1}^\vartriangle)} ,~ e_{(2E_{1,2}^\vartriangle)}$$ are
$1, 3, 0, 2$ respectively, and the degrees of $x_1, x_2, x_2^{-1}
$ are $2, 4, -4$ respectively. Since $x_1$ and $x_2$ both have
even degrees, we only need to prove that
$e_{(2E_{1,1}^\vartriangle)} ,~ e_{(2E_{1,2}^\vartriangle)}$ are
linearly independent and
$e_{(E_{1,1}^\vartriangle+E_{1,2}^\vartriangle)},~e_{(E_{1,2}^\vartriangle+E_{1,3}^\vartriangle)}~$
are linearly independent.

Suppose there is an equation $$f\cdot
e_{(2E_{1,1}^\vartriangle)}+g\cdot e_{(2E_{1,2}^\vartriangle)}=0,
$$ where $f, g\in \bbq[x_1, x_2, x_2^{-1}]$. We want to show that $f=g=0$.
 Since $x_1\cdot
e_{(2E_{1,1}^\vartriangle)}=x_1$ and $x_2\cdot
e_{(2E_{1,1}^\vartriangle)}=x_2$, then $f\cdot
e_{(2E_{1,1}^\vartriangle)}=f$ and
\begin{equation}\label{li}f+g\cdot e_{(2E_{1,2}^\vartriangle)}=0. \end{equation} Suppose
$f=\sum_Aa_Ae_A$, where $a_A\in\bbq$ and $A\in\cnr$. Then by
\eqref{ee} each term $e_A$ in $f$ satisfies
$\col(A)=\col(2E_{1,1}^\vartriangle)=\lb$. Suppose $g\cdot
e_{(2E_{1,2}^\vartriangle)}=\sum_{A'}b_{A'}e_{A'}$, where
$b_{A'}\in\bbq$ and $A'\in\cnr$. Then by \eqref{ee} each term
$e_{A'}$ in $g\cdot e_{(2E_{1,2}^\vartriangle)}$ has the property
that $\col(A')=\col(2E_{1,2}^\vartriangle)=\mu$. Since the set
$\{e_A~|~A\in\cnr\}$ is a basis for $\s$, the basis corresponding
to matrices of different columns are linearly independent. This
implies that $f=0$ and $g\cdot e_{(2E_{1,2}^\vartriangle)}=0$.
Then $b_{A'}=0$ for each $A'$ occurring in $g\cdot
e_{(2E_{1,2}^\vartriangle)}$. By proposition \eqref{mu11}, we get
that $g=\sum_{A'[1]}b_{A'}e_{A'[1]}$, where $A'[1]$ is the matrix
that $A'[1]$'s $i$th-column is the same as that $A'$'s
$(i+1)$th-column. This proves that $g=0$ and
$e_{(2E_{1,1}^\vartriangle)} ,~ e_{(2E_{1,2}^\vartriangle)}$ are
linearly independent.

Suppose there is an equation \begin{equation}\label{dd}f'\cdot
e_{(E_{1,1}^\vartriangle+E_{1,2}^\vartriangle)}+g'\cdot
e_{(E_{1,2}^\vartriangle+E_{1,3}^\vartriangle)}=0,
\end{equation} where $f', g'\in \bbq[x_1,
x_2, x_2^{-1}]$ and at least one of $f', g'$ are not zero. We want
to find a contradiction. Since $x_2$ is invertible, we assume each
term containing $x_2$ in $f'$ and $g'$ has positive degrees in
$x_2$.

Now suppose $f'=\sum_{a,b}m_{a,b} x_1^ax_2^b$ and
$g'=\sum_{c,d}n_{c,d} x_1^cx_2^d$, where $m_{a,b}, ~n_{c,d}\in
\bbq$. By Proposition \eqref{mu3} and Proposition \eqref{mu4}, we
get
\begin{equation}\label{d}x_1^ax_2^b\cdot
e_{(E_{1,1}^\vartriangle+E_{1,2}^\vartriangle)}=\sum_{0\leq k\leq
2a}p_ke_{(E_{1,2b+1+k}^\vartriangle+E_{1,2b+2+2a-k}^\vartriangle)}
, \end{equation}
\begin{equation}\label{d}x_1^cx_2^d\cdot
e_{(E_{1,2}^\vartriangle+E_{1,3}^\vartriangle)}=\sum_{0\leq l\leq
2c}h_le_{(E_{1,2d+2+l}^\vartriangle+E_{1,2d+3+2c-l}^\vartriangle)}
,\end{equation} where \begin{equation}\label{z}\begin{cases}
p_k=1,  &k=0,\;\text{or}\; 2a ;\\
p_k>1,& 1\leq k\leq 2a-1;\\
\end{cases}\;\;
 \begin{cases}
h_l=1,  &l=0,\;\text{or}\; 2c ;\\
h_l>1,& 1\leq l\leq 2c-1.\\
\end{cases}\end{equation}

Suppose $f'\neq 0$. Let $m_{a_0,b_0}x_1^{a_0}x_2^{b_0}$ be any
nonzero term appearing in $f'$ with $b_0$ minimal. Suppose
$e_{(E_{1,2b_0+1}^\vartriangle+E_{1,2b_0+2+2a_0}^\vartriangle)}$
appears in $f'\cdot
e_{(E_{1,1}^\vartriangle+E_{1,2}^\vartriangle)}$ with coefficient
$s$, then by \eqref{z}, $s\geq m_{a_0,b_0}\neq 0$. By \eqref{dd},
there is some $d_0 $ and some $l$ with $0\leq l\leq 2c$ such that
$$e_{(E_{1,2b_0+1}^\vartriangle+E_{1,2b_0+2+2a_0}^\vartriangle)}=
e_{(E_{1,2d_0+2+l}^\vartriangle+E_{1,2d_0+3+2c-l}^\vartriangle)}.$$
 Then we get that $e_{(E_{1,2d_0+2}^\vartriangle+E_{1,2d_0+3+2c}^\vartriangle)}$
appears in  $g'\cdot
e_{(E_{1,2}^\vartriangle+E_{1,3}^\vartriangle)}$ with nonzero
coefficient by \eqref{z}. By $b_0$ is minimal, we cannot find this
term in $f'\cdot e_{(E_{1,1}^\vartriangle+E_{1,2}^\vartriangle)}$.
 This
contradicts \eqref{dd}.

Similarly we can find a contradiction if we suppose $g'\neq 0$.
This proves that
$e_{(E_{1,1}^\vartriangle+E_{1,2}^\vartriangle)},~e_{(E_{1,2}^\vartriangle+E_{1,3}^\vartriangle)}~$
are linearly independent.

Then $\Pi_1$ is a basis of $e_\lb \s$ over $e_\lb \s e_\lb$. In a
similar way we can prove that the set
$$\Pi_2=\{e_{(E_{1,1}^\vartriangle+E_{2,1}^\vartriangle)},~
e_{(E_{2,1}^\vartriangle+E_{3,1}^\vartriangle)},~
e_{(2E_{1,1}^\vartriangle)} ,~ e_{(2E_{2,1}^\vartriangle)} \}$$
forms a basis of $\s e_\lb$ as a free $e_\lb \s e_\lb$-module.

\end{pf}

\begin{prop}\label{bimodule} There is an $\s$-bimodule isomorphism $$\alpha:
\s e_\lb\otimes_{B}e_\lb \s\longrightarrow \s e_\lb \s,$$ and
$J=\s e_\lb \s$ is projective as left and right $\s$-module.
\end{prop}

\begin{pf}

 Let $B=e_\lb \s e_\lb$ and let $\{b_i\}_{i\in I} $ be a $\bbq$-basis of
$B$. There is a canonical epimorphism $$\alpha:~~\s
e_\lb\otimes_{B}e_\lb \s\longrightarrow \s e_\lb \s$$
$$e_Ab_i\otimes b_{i'} e_{C} \longrightarrow e_Ab_i b_{i'} e_{C},~~
e_A\in \Pi_2, e_C\in \Pi_1, i,i'\in I.$$

The set $$\Omega=\{e_A b_i\otimes e_C~|~ e_A\in \Pi_2, e_C\in
\Pi_1, i \in I\}$$ is a basis of $\s e_\lb\otimes_{B}e_\lb \s$
over $\bbq$. Let
$$\Omega'=\{e_A b_i e_C~|~ e_A\in \Pi_2, e_C\in \Pi_1, i\in I\}.$$

We prove $\alpha$ is injective by showing that elements in
$\Omega'$ are $\bbq$-linearly independent.

Recall that
$$\Pi_1=\{e_{(E_{1,1}^\vartriangle+E_{1,2}^\vartriangle)},~
e_{(E_{1,2}^\vartriangle+E_{1,3}^\vartriangle)},~
e_{(2E_{1,1}^\vartriangle)} ,~ e_{(2E_{1,2}^\vartriangle)} \}$$
and
$$\Pi_2=\{e_{(E_{1,1}^\vartriangle+E_{2,1}^\vartriangle)},~
e_{(E_{2,1}^\vartriangle+E_{3,1}^\vartriangle)},~
e_{(2E_{1,1}^\vartriangle)} ,~ e_{(2E_{2,1}^\vartriangle)} \}.$$

Let $\pi_l$ for $1\leq l\leq 4$ denote the elements in $\Pi_1$ and
let $\pi_l'$ for $1\leq l\leq 4$ denote the elements in $\Pi_2$.
Define $$\Omega'_{l,m}=\{\pi_l' b_i \pi_m~|~ i\in I\},$$ then
$$\Omega'=\cup_{1\leq l,m\leq 4}\Omega'_{l,m}.$$

The elements in $\Omega'_{l,m}$ are $\bbq$ combinations of basis
$e_A$ such that $$\col(A)=\col(C), \row(A)=\row(D) \;\text{for}\;
1\leq l,m\leq 4,$$ where $\pi_l'=e_C$ and $\pi_m=e_D$. Since
$\{e_A~|~A\in\Theta_\vartriangle(2,2)\}$ is a basis of $\s$, we
can divide $\Omega'$ into disjoint union of subsets according to
the row and column vectors of the matrices corresponding to the
basis elements, such that elements in different subsets are
linearly independent. It suffices to prove elements in the
following sets are linearly independent.
$$ \Omega'_{l,m}, 3\leq l,m\leq 4,$$
$$\bigcup_{1\leq l,m\leq 2}\Omega'_{l,m},$$
$$ \Omega'_{3,1}\cup\Omega'_{3,2},$$
$$ \Omega'_{4,1}\cup\Omega'_{4,2},$$
$$ \Omega'_{1,3}\cup\Omega'_{2,3},$$
$$ \Omega'_{1,4}\cup\Omega'_{2,4},$$

 We only prove
 elements in $ \Omega'_{3,1}\cup\Omega'_{3,2}$ are linearly independent. Other cases can
be proved similarly. Suppose
$$e_{(2E_{1,1}^\vartriangle)}\cdot f\cdot e_{(E_{1,1}^\vartriangle+E_{1,2}^\vartriangle)}+
e_{(2E_{1,1}^\vartriangle)}\cdot g\cdot
e_{(E_{1,2}^\vartriangle+E_{1,3}^\vartriangle)}=0,$$ where $f, g$
are $\bbq$-linear combinations of elements in $\{b_i\}_{i\in I}$.
Then $$e_{(2E_{1,1}^\vartriangle)}\cdot\biggl( f\cdot
e_{(E_{1,1}^\vartriangle+E_{1,2}^\vartriangle)}+
 g\cdot
e_{(E_{1,2}^\vartriangle+E_{1,3}^\vartriangle)}\biggr)=0.$$ by
Proposition  \eqref{mu2},
$$ f\cdot
e_{(E_{1,1}^\vartriangle+E_{1,2}^\vartriangle)}+
 g\cdot
e_{(E_{1,2}^\vartriangle+E_{1,3}^\vartriangle)}=0.$$ Then we can
get $f=g=0$, since $
e_{(E_{1,1}^\vartriangle+E_{1,2}^\vartriangle)},
e_{(E_{1,2}^\vartriangle+E_{1,3}^\vartriangle)}$ belong to
$\Pi_1$.

This proves that $J\cong\s e_\lb\otimes_{B}e_\lb \s$ as
$\s$-bimodule. Now we prove $\s e_\lb\otimes_{B}e_\lb \s$ is
projective as left and right $\s$-module.

 Since $ e_\lb\s$ and $\s e_\lb$ are free modules over
$B$, $\s e_\lb\otimes_{B}e_\lb \s$ is projective as $\s$-bimodule.
Then $J$ is projective as $\s$-bimodule follows from the
isomorphism constructed above.

\end{pf}

\begin{thm}\label{af} $\s$ is affine cellular over $\bbq$.
\end{thm}

\begin{pf}
Let $\tau$ be the anti-involution of $\s$ given by $\tau(e_A)
=e_{A'}$, where $A\in \Theta_\vartriangle(2,2)$ and $A'$ denotes
the transpose of $A$. Let $J=\s e_\lb \s$. It is obvious that
$\tau(J)=J$. We first show that $J$ is an affine cell ideal of
$\s$.

 Let $$\Pi_1=\{e_{(E_{1,1}^\vartriangle+E_{1,2}^\vartriangle)},~
e_{(E_{1,2}^\vartriangle+E_{1,3}^\vartriangle)},~
e_{(2E_{1,1}^\vartriangle)} ,~ e_{(2E_{1,2}^\vartriangle)} \}$$
 and $$\Pi_2=\{e_{(E_{1,1}^\vartriangle+E_{2,1}^\vartriangle)},~
e_{(E_{2,1}^\vartriangle+E_{3,1}^\vartriangle)},~
e_{(2E_{1,1}^\vartriangle)} ,~ e_{(2E_{2,1}^\vartriangle)} \}$$
 be the same as in Lemma \eqref{main}.
It is easy to find that
\begin{equation}\label{dual}\Pi_1=\{\tau(e_{A})~|~ e_A\in
\Pi_2\}.\end{equation} Let $V$ be the free $\bbq$-module on the
basis $\Pi_2$. Suppose
$$ ~~\triangle=V\otimes_\bbq B, ~~\triangle'=B\otimes_\bbq V,$$
where $B= e_\lb\s e_\lb$. Then by Lemma \eqref{main} we have the
following isomorphism of $\s$-$B$-bimodules
$$\varphi:~~\triangle\longrightarrow\s e_\lb ,$$
$$ \varphi: v\otimes b\longrightarrow
v\cdot b,$$ where $v\in V, b\in B $.
 And by \eqref{dual} we
have the following isomorphism of $B$-$\s$-bimodules
$$\psi:~~\triangle'\longrightarrow e_\lb \s,$$
$$ \psi:b\otimes v\longrightarrow b\cdot\tau(v),$$
where $ b\in B, v\in V $.

 Then by Lemma \eqref{main} we have the
following $\s$-bimodule isomorphism:
\begin{equation}\label{biiso}\alpha: \triangle \otimes_B \triangle'\longrightarrow J,\end{equation}
$$\alpha(v_1\otimes b_1\otimes
b_2\otimes v_2)=v_1\cdot b_1\cdot b_2\cdot\tau(v_2),$$ where $b_1,
b_2\in B$ and $v_1, v_2\in V$.

Now define $$f:\triangle \otimes_B \triangle'\longrightarrow
\triangle \otimes_B \triangle',$$ $$f: v_1\otimes b_1\otimes
b_2\otimes v_2\longrightarrow v_2\otimes \tau(b_2)\otimes
\tau(b_1)\otimes v_1,$$ where $b_1, b_2\in B$ and $v_1, v_2\in V$.

Then we have the following commutative diagram:
$$\xymatrix{ J \ar[d]^-{\tau}& \triangle \otimes_B \triangle'\ar[l]^-{\alpha}\ar[d]^-{
f}\\
J &\triangle \otimes_B \triangle' \ar[l]^-{\alpha}}$$

Let $\alpha^{-1}$ be the inverse of the bimodule isomorphism
$\alpha$ given in \eqref{biiso}. Since $\tau \alpha=\alpha f$, we
get that $\alpha^{-1}\tau=f \alpha^{-1}$, i.e. we have the
following commutative diagram:

$$\xymatrix{ J \ar[d]^-{\tau}\ar[r]^-{\alpha^{-1}}& \triangle \otimes_B \triangle'\ar[d]^-{f}\\
J\ar[r]^-{\alpha^{-1}} &\triangle \otimes_B \triangle' }$$

  This shows that $J$ is an affine cell ideal of $\s$.
 By Lemma
\eqref{maa}, $\s/J\cong \bbq[x,x^{-1}]$. Then  $\s/J$ is an affine
cellular algebra over $\bbq$. Now we have a chain of ideals
$$0\subset J\subset \s$$ and a decomposition of $\s$ as a $\bbq$-vector space( In fact,
this is an decomposition of $\s$-module):
$$\s=J\bigoplus \s/J.$$

By Lemma \eqref{maa}, $\tau(\s/J)=\s/J$. This implies that this
chain and decomposition satisfy the conditions given in
\eqref{def-aff}. This proves that $\s$ is an affine cellular
algebras over $\bbq$.

\end{pf}

\begin{cor}\label{pr} $\s$ is of finite global
dimension over $\bbq$.
\end{cor}

\begin{pf}

By Lemma 4.5 in \cite{kx}, if $R$ is a ring and $e=e^2\in R$, if
$ReR$ is projective as left $R$-module, then the global dimension
of $R$ is finite if and only if the global dimension of $(R/ReR)$
and the global dimension of $(eRe)$ are finite.

By Proposition \eqref{bimodule}, $J=\s e_\lb\s $ is projective as
left $\s$-module. Since $\s/J\cong \bbq[x, x^{-1}]$ and $e_\lb \s
e_\lb\cong \bbq[x_1, x_2, x_2^{-1}]$ by Lemma \eqref{maa} and
\eqref{m}. Then $\s/J$ and $e_\lb \s e_\lb$ are both of finite
global dimensions. This proves that $\s$ is of finite global
dimension over $\bbq$.

\end{pf}

\bigskip

\noindent\small {\bf Acknowledgements}\;\; The author is grateful
to Professor Bangming Deng for his guidance and stimulating
discussions.

\end{document}